\documentclass[12pt]{amsart}
\usepackage{url,hyperref,lineno,microtype,subcaption}
\usepackage{mathrsfs}
\usepackage{graphicx}
\usepackage{amsmath}
\usepackage{amsfonts}
\usepackage{amsthm}
\usepackage{amscd}
\usepackage{amssymb}
\usepackage{color}
\usepackage{enumerate}
\usepackage{tikz}
\usetikzlibrary{calc}
\usetikzlibrary{shapes.geometric}
\usetikzlibrary{lindenmayersystems}
\usepackage{tikzit}
\tikzstyle{white node}=[fill=white, draw=black, shape=circle]
\tikzstyle{purple}=[fill={rgb,255: red,117; green,27; blue,196}, draw=black, shape=circle]
\tikzstyle{green}=[fill={rgb,255: red,0; green,197; blue,0}, draw=black, shape=circle]
\tikzstyle{orange}=[fill={rgb,255: red,255; green,128; blue,0}, draw=black, shape=circle]
\tikzstyle{black}=[fill=black, draw=black, shape=circle]

\tikzstyle{black edge}=[-, draw=black, fill=black]
\tikzstyle{dark grey}=[-, fill={rgb,255: red,128; green,128; blue,128}, draw={rgb,255: red,128; green,128; blue,128}]
\tikzstyle{red edge}=[-, fill=none, draw=red]
\tikzstyle{blue edge}=[-, draw=blue]
\tikzstyle{green edge}=[-, draw=green]
\tikzstyle{yellow edge}=[-, draw={rgb,255: red,222; green,222; blue,0}]

\newtheorem{theorem}{Theorem}[section]
\newtheorem{lemma}[theorem]{Lemma}

\newtheorem{note}[theorem]{Note}

\newtheorem{prop}[theorem]{Proposition}
\newtheorem{observation}[theorem]{Observation}

\newtheorem{corollary}[theorem]{Corollary}

\theoremstyle{definition}
\newtheorem{definition}[theorem]{Definition}
\theoremstyle{remark}
\newtheorem{remark}[theorem]{Remark}

\def\ec{{\rm ec}\hskip0.02cm}

\newcommand{\str}{{\tt Stretch}\hskip0.02cm}
\newcommand\remove[1]{}

\begin{document}

\title{Quantitative characteristics of cycles and their relations with stretch and spanning tree congestion}

\author{Florin Catrina
\and
 Rainah Khan
 \and
  Isaac Moorman
  \and
  Mikhail Ostrovskii
  \and
   Lakshmi Iswara Chandra Vidyasagar   }

\address{Department of Mathematics and  Computer Science \\
St. John's University\\
Queens, New York 11439}

\email{catrinaf@stjohns.edu; rainah.khan19@my.stjohns.edu; isaac.moorman18@my.stjohns.edu; ostrovsm@stjohns.edu; iswaracl@stjohns.edu}

\thanks{
\noindent Research of Mikhail Ostrovskii and of two undergraduate
students participating in this project - Rainah Khan and Isaac Moorman - was partially supported by  NSF DMS-1700176 and DMS-1953773.\\
The authors gratefully acknowledge the free hosting and support of
the group's discussion board on the
\href{https://reu-coauthor.csail.mit.edu/}{REU-Coauthor}
platform.}



\begin{abstract}
The main goal of this article is to introduce new quantitative
characteristics of cycles in finite simple connected graphs and to
establish relations of these characteristics with the {\em
stretch} and {\em spanning tree congestion} of graphs. The main
new parameter is named the {\em support number}. We give a
polynomial approximation algorithm for the support number with the
aid of yet another characteristic we introduce, named the {\em
cycle width} of the graph.
\end{abstract}

\date{\today}

\maketitle

%This is version 1

\section{Introduction}

This paper is devoted to connections between several
characteristics of finite connected simple graphs. Two such
characteristics are the well-known notion of {\em stretch}
(introduced in \cite{PU89}, see surveys in \cite{CC95, FK01, LW08,
Pel00}) and that of {\em spanning tree congestion} (introduced in
\cite{Ost04}, see also \cite{LO10, Ost10, Ota11, Ota20}). These
two notions are closely related for plane graphs (see Lemma
\ref{L:CongStr} below).

It is known that the problem of finding the minimum stretch
spanning tree in a graph $G$ is NP-hard \cite{CC95,FK01}.
An approximation algorithm with approximations ratio
$O(\log n)$ was found for it in \cite{EP09}.

The main goal of this paper is to introduce and to study a new
related characteristic, namely that of  {\it support number.} We
establish relations between this characteristic, stretch, and
spanning tree congestion. The main result of the paper is an
approximation algorithm with constant approximation ratio for
computing the support number. This is done by introducing yet
another characteristic {\em the cycle width} of a graph (denoted
$W(G)$) which is computable in polynomial time and is equivalent
to the support number $k$ in the sense
\[ W(G)/3 \le k \le W(G)+4.\]

We conjecture that finding the support number of $G$ is an NP-hard
problem.

\section{Congestion and Stretch}\label{S:cs}

Our Graph Theory terminology and notations are standard (see
\cite{BM08,CH91}). Let $G$ be a finite, simple, connected graph. We
denote by $V(G)$ the vertex set of $G$, and by $E(G)$ the edge set
of $G$. Recall that a {\em spanning subgraph} of $G$ is a subgraph
that contains all the vertices $V(G)$ of the original graph. A
{\em spanning tree} of $G$ is a connected, spanning subgraph,
which does not have any cycles.

We use the following definitions introduced in \cite{Ost04} (some
of them were known before, see \cite{Sim87}).
 Let $T$ be a spanning tree in $G$.
\begin{definition}[\cite{Ost04}]
\label{D:ed_con_givenT}
 For an edge $e$ of $T$
let $A_e$ and $B_e$ be the vertex sets of the two components of $T-e$.
Denote by $e_G(A_e, B_e)$ the number of edges in $G$ that connect a vertex in $A_e$ to a vertex in $B_e$.

\begin{figure}[h!]
\centering
\begin{subfigure}{.5\textwidth}
  \centering
  \includegraphics[width=.6\linewidth]{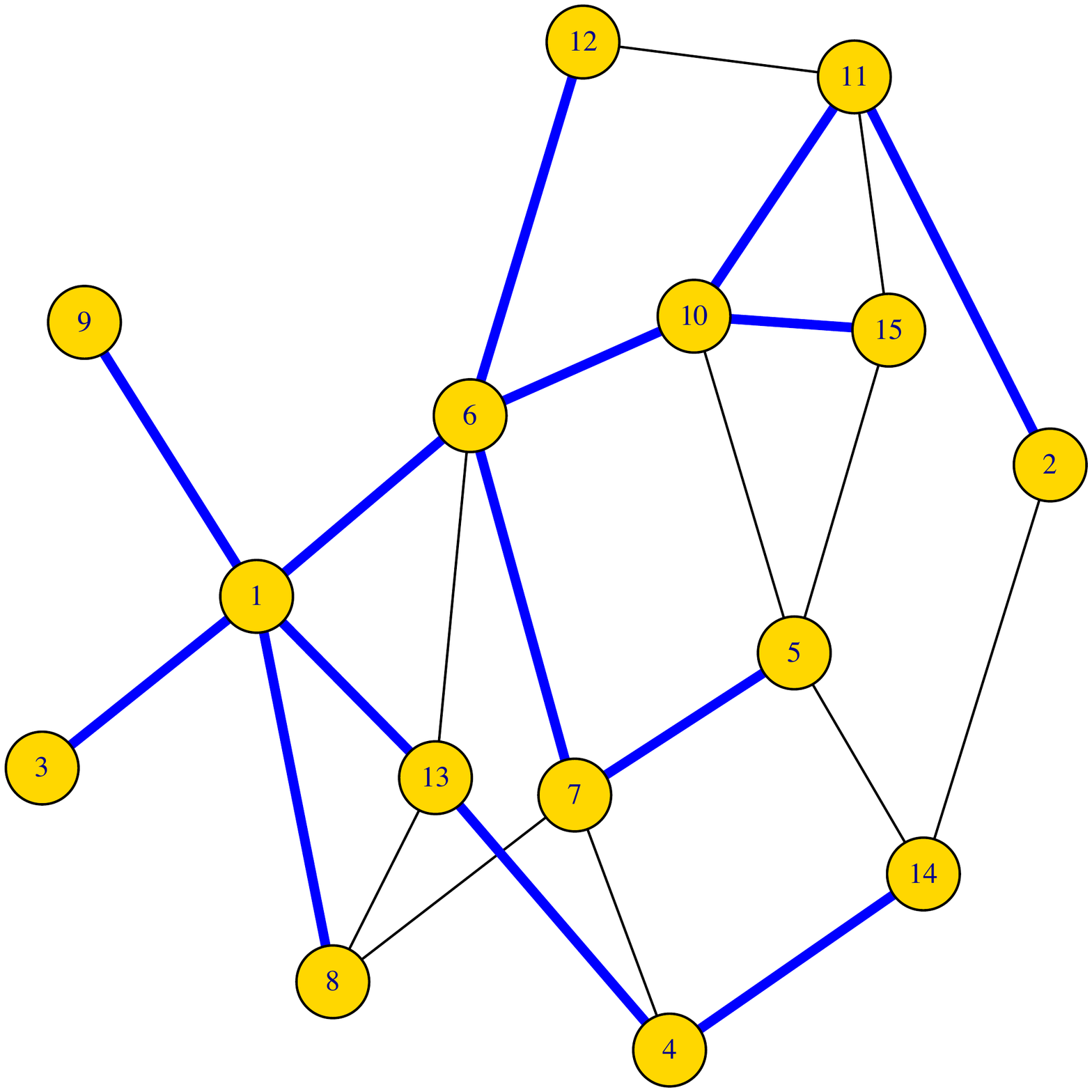}
  \caption{Spanning tree consists of blue edges}
  %\label{fig:sub1}
\end{subfigure}%
\begin{subfigure}{.5\textwidth}
  \centering
  \includegraphics[width=.6\linewidth]{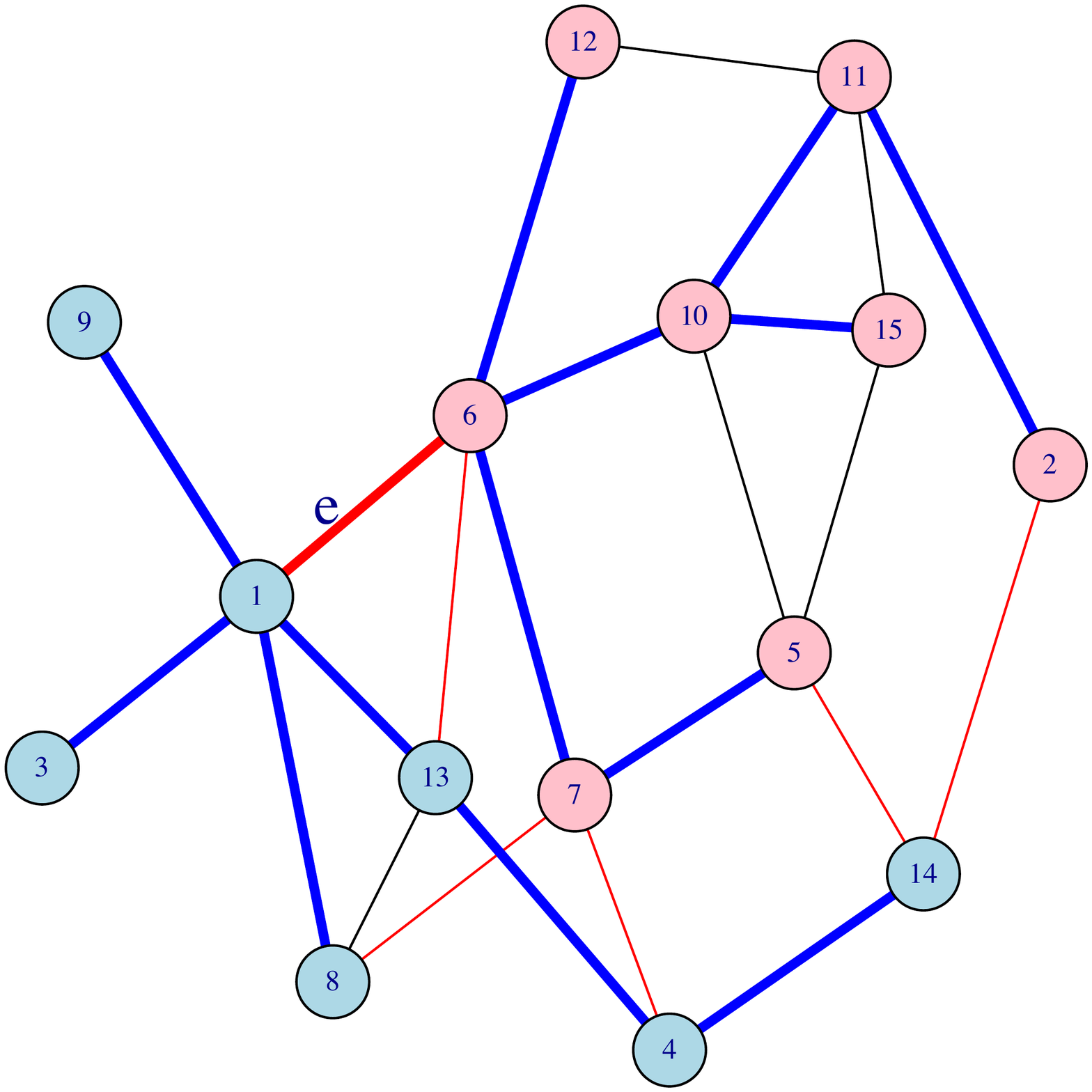}
  \caption{$ A_e=$ blue vertices,  $ B_e= $ pink vertices.}
 % \label{fig:sub2}
\end{subfigure}
\caption{Red edges  show $e_G(A_e, B_e) = 6$. }
\label{fig:congest}
\end{figure}

\noindent
The {\it edge congestion} of $G$ in $T$
is
$$
\ec(G:T)=\max_{e\in E(T)} e_G(A_e,B_e).
$$
The {\it spanning tree congestion} of $G$
is
$$
s(G)=\min_{T} \ec(G:T).
$$
\end{definition}
In short, the {\em spanning tree congestion} of $G$ is the min-max
characteristic
\[ s(G)=\min_{T}
\max_{e\in E(T)} e_G(A_e,B_e).
\]

\begin{definition}[\cite{PU89}]
 If $H$ is a connected spanning subgraph in $G$, then the {\it
stretch of $G$ in $H$} is defined by
%\begin{equation}\label{E:str}
\[ \str(G:H)=\max_{u,v\in V(G)}\frac{d_H(u,v)}{d_G(u,v)}\]
%\end{equation}
\end{definition}
It suffices to test the stretch only for adjacent vertices $u,v$
in $G$. That is,
\[\str(G:H) =\max_{uv\in E(G)} d_H(u,v). \]
The {\em stretch} of $G$ is defined by minimizing with respect to the spanning trees $T$ in $G$
\[ \sigma(G) = \min_T \str(G:T).\]
In short, the {\em stretch} of $G$ is the min-max
characteristic
\[ \sigma(G)=\min_{T}
\max_{uv\in E(G)} d_T(u,v).
\]
\begin{definition}
The {\em dual graph} $G^*$ of a plane graph $G$ is the multigraph
whose vertices correspond to the faces of $G$, including the
exterior face. Every edge $e^*$  joining two vertices of $G^*$
corresponds to an edge $e$ of $G$ in the common boundary of
the corresponding faces of $G$, and vice versa.
 \end{definition}

If two faces have several common edges in their boundaries, the
corresponding edges are multiple edges in $G^*$.

\begin{definition}
If $T$ is a spanning tree of a plane graph $G$, then the {\em dual
spanning tree} $T^\sharp$ is defined as the spanning subgraph of
$G^*$ such that $e^*\in E(T^\sharp)$ if and only if $e\notin
E(T)$.
\end{definition}

The connection between stretch and spanning tree congestion is revealed in the case of planar graphs
by the following lemma (Lemma 8 in \cite{LLO14}).

\begin{lemma}
\label{L:CongStr} Let $G$ be a connected planar graph.
\begin{itemize}
\item[(a)] If $T$ is a spanning tree in $G$ and $T^\sharp$
is its dual spanning tree, then
$$\ec(G:T)=\str(G^*:T^\sharp)+1.$$
\item[(b)] $s(G)=\sigma(G^*)+1$.
\end{itemize}
\end{lemma}

\begin{note} It is proved in \cite{FK01} that
determination of the least $t$ for which a planar graph has a
spanning tree $T$ with $\str(G:T)=t$ is NP-hard. Combining this
with Lemma \ref{L:CongStr} we get that the problem of computing
$s(G)$ for planar graphs is also NP-hard. This fact was also
observed in \cite{BFGOV12,Low10}. The paper \cite{BFGOV12}
contains many interesting results on the complexity of the
spanning tree congestion.
\end{note}

\section{Statements of the main results}

Since there are no efficient algorithms for calculating stretch
and congestion for general graphs, it is highly desirable to find
polynomially computable parameters which can be used to estimate
the two characteristics. This is the main motivation for our
study.

%%%%%%

Presence of a cycle in a graph implies that $\sigma(G)\ge 2$.
Without additional restrictions on the cycle, it does not imply
more; complete graph has stretch $2$ and contains cycles of all
possible lengths. To get stronger estimates for stretch from
below, we need cycles which, in some coarse metric sense, resemble
circles rather than trees. The strongest notion of this type is
that of a ``long'' isometric cycle in a graph.

\begin{definition} Let $C$ be a cycle in a graph $G$.
We say that $C$ is an {\it isometric cycle} in $G$ if for any two
vertices  $u,v$  of $C$,
\[ d_{C}(u,v)=d_G(u,v).\]
\end{definition}

We will be interested also in the following weakening of the
notion of isometric cycle:

\begin{definition}
Let $\alpha\in(0,1)$. We say that $C$ is a $(\alpha,1)$-{\it
bilipschitz cycle} in $G$ if
\[\forall u,v\in C\quad \alpha d_{C}(u,v)\le d_G(u,v)\le d_{C}(u,v).\]
\end{definition}

Isometric and bilipschitz classes of cycles are proper subclasses
of the class of supported cycles introduced in the next
definition. The class of $k$-supported cycles seems to be the most
relevant in the study of the stretch of graphs.

\begin{definition}\label{D:kSuppCyc} If a cycle $C$ in a graph $G$  can
be partitioned into three edge disjoint paths $I_1,I_2$, and $I_3$
with  $I_1\cap I_2$, $I_2\cap I_3$,  and $I_3\cap I_1$ containing
one vertex each, and in such a way that for every triple
$(u_1,u_2,u_3)$ of vertices satisfying $u_i\in I_i$, $i=1,2,3$, we
have
\[\max_{i,j\in\{1,2,3\}} d_G(u_i,u_j)\ge k,\]
we say that $C$ is a \emph{$k$-supported cycle} contained in $G$.
\end{definition}

\begin{figure}[h!]
\includegraphics[width=.5\textwidth, angle=90]{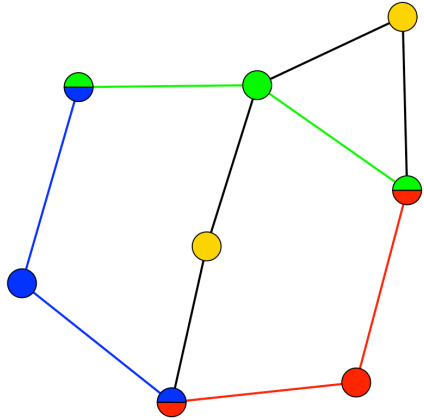}
  \caption{The cycle  obtained by concatenating the green, red, and blue paths is $2$-supported in $G$.}
    \label{fig:kSuppCyc}
\end{figure}

\begin{definition}\label{D:SuppNum}
For a finite, simple, connected graph $G$, define its {\em support
number} as the largest integer $k$ for which $G$ contains a
$k$-supported cycle. If $G$ is a tree,  the support number is
defined to be zero.
\end{definition}

\begin{remark}
If $G$ is not a tree, then any of its cycles is $1$-supported.
This can be seen by considering  any two consecutive edges of $C$ as $I_1$ and $I_2$, while the remaining part of the cycle is
the path $I_3$.
Therefore the support number  of any non-tree is some $k\ge 1$.
\end{remark}

The following theorem states that the existence of a $k$-supported
cycle in a graph $G$ leads to an estimate from below for the
stretch of $G$.
\begin{theorem}\label{T:kSuppImplStr} If a graph $G$ contains a $k$-supported cycle,
then $\sigma(G)\ge k$.
\end{theorem}
If the graph contains large isometrically embedded cycles, the
result above implies the following estimate from below for the
stretch.
\begin{corollary}
\label{C:iso_embd} If a graph $G$ contains an isometrically
embedded cycle $C$ of length $n$, then $\sigma(G)\ge
\left\lceil\frac{n}3\right\rceil$.
\end{corollary}

\begin{remark}
Similarly, if a graph $G$ contains an $(\alpha,1)$-bilipschitz
embedded cycle $C$, then $C$ is $k$-supported with $k=
 \alpha\left\lceil\frac{n}3\right\rceil$ and hence $\sigma(G)\ge
 \alpha\left\lceil\frac{n}3\right\rceil$.
\end{remark}

We find also an approximate algorithm for computing the support
number. Namely, for a graph $G$ we describe below
(Definition~\ref{D:cw}) a polynomially computable numerical
parameter $W(G)$  and prove two theorems.

\begin{theorem}\label{T:(A)} Each graph $G$ contains a
$W(G)/3$-supported cycle.
\end{theorem}

\begin{theorem}\label{T:(B)}
Each graph $G$ containing a $k$-supported cycle satisfies $W(G)\ge
k-4$.
\end{theorem}

\section{$k$-supported cycles}
\label{S:k_sp_Cy}

Before proving Theorem~\ref{T:kSuppImplStr}, we estimate the
support number for rectangular and triangular grids.
\begin{prop}
\label{P:rectangular}
 For the rectangular grid $G= P_m\times P_n$ with
$2\leq m\leq n$, the exterior rectangle is an $(m-1)$-supported cycle.
\end{prop}
\proof
We color the left side (this is path $I_1$ with $(m-1)$ vertical edges) in red, the top side (this is path $I_2$ with $(n-1)$ horizontal edges) in green,
and the remaining two sides (this is path $I_3$) in blue.

\begin{figure}[h!]
\includegraphics[width=.4\textwidth, angle=270]{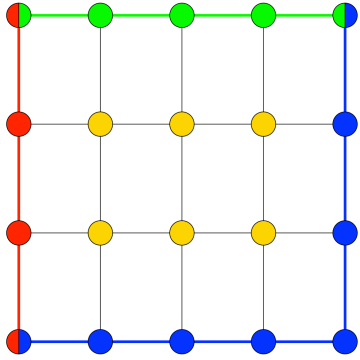}
  \caption{Rectangular grid $P_4\times P_5$.}
    \label{fig:re_grid}
\end{figure}

The blue vertex $u_3$ is either on the bottom horizontal path and therefore at distance at least $(m-1)$ to the green vertex $u_2$,
or on the right vertical path and  therefore at distance at least $(n-1)$ to the red vertex $u_1$.
In either case,
\[\max_{i,j\in\{1,2,3\}} d_G(u_i,u_j)\ge (m-1).\]
\endproof

\begin{prop}
\label{P:triangular}
 For the triangular grid $G= T_{n}$ with
$n\ge 2$, the exterior triangle is an $ \lfloor
\frac{n-1}{2}\rfloor$-supported cycle.
\end{prop}
\proof We prove in the Euclidean case that an equilateral triangle
of side $\ell$  has the property that for any three points
$u_1,u_2,u_3$, one on each side, it holds that
\begin{equation}
\label{E:l/2}
\max_{i,j\in\{1,2,3\}} d(u_i,u_j)\ge \frac \ell 2.
\end{equation}
Here, $d(u_i,u_j)$ is the Euclidean distance in the plane. Since the
graph distance on the triangular grid dominates the Euclidean
distance, inequality \eqref{E:l/2} implies that for any three
vertices  $u_1,u_2,u_3$   of the grid, one on each side, it holds
that
\[\max_{i,j\in\{1,2,3\}} d_{G}(u_i,u_j)\ge  \frac{n-1}{2}.\]

Therefore, let us consider the Euclidean case of an equilateral triangle of side $\ell$.
Assume by contradiction that there exist  $u_1,u_2,u_3$,  one on each side, such that
\begin{equation}
\label{E:contr}
\max_{i,j\in\{1,2,3\}} d(u_i,u_j)< \frac \ell 2.
\end{equation}

\begin{figure}[h!]
\centering
\begin{subfigure}{.5\textwidth}
  \centering
  \includegraphics[width=.8\linewidth]{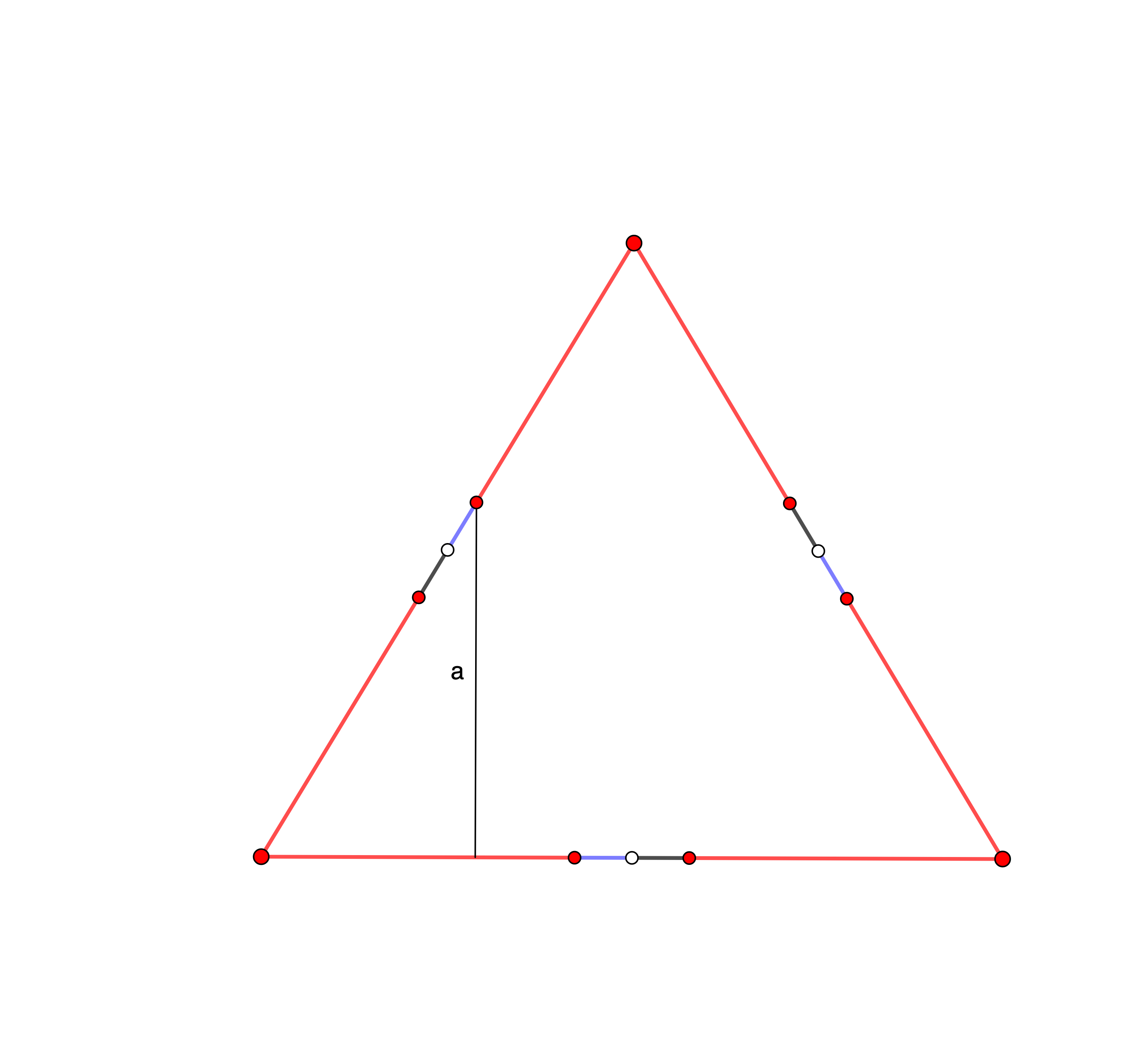}
  \caption{Non-admissible region in red}
  %\label{fig:sub1}
\end{subfigure}%
\begin{subfigure}{.5\textwidth}
  \centering
  \includegraphics[width=.9\linewidth]{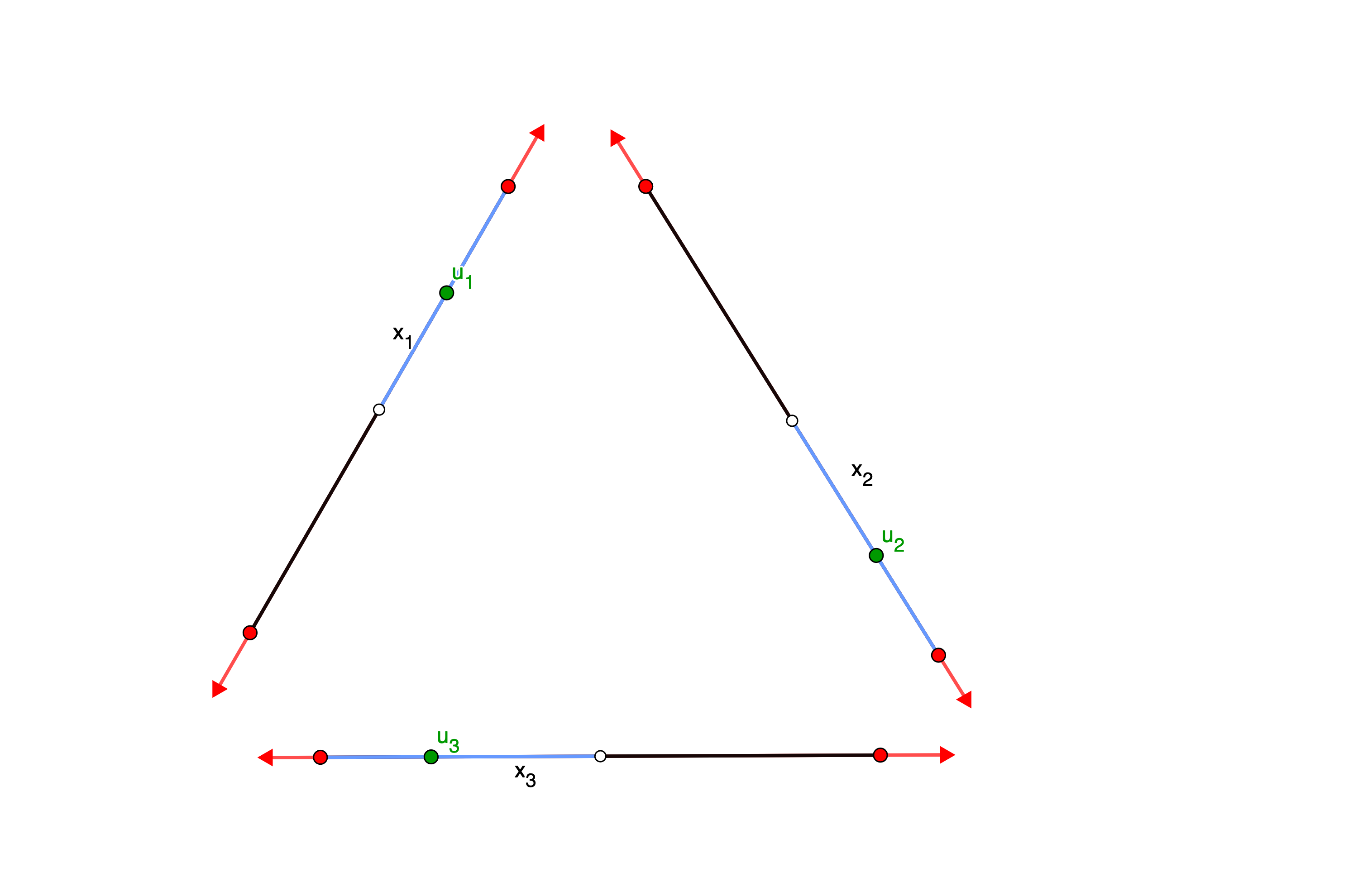}
  \caption{Admissible region}
 % \label{fig:sub2}
\end{subfigure}
\caption{Equilateral triangle  of side $\ell=2a$}
\label{F:equilat}
\end{figure}

Clearly, none of the points can be in the red region because the distance from a red wedge to the opposite side is $\ell/2$.

Also, none of the points can be a  midpoint on its side,
for if for example $u_1$ is midpoint of its side, then $u_2$ has to be in the black region and $u_3$ in the blue,
therefore $d(u_2,u_3)\ge \ell/2$.

We claim that  $u_1,u_2,u_3$ are all either in the blue intervals or all in the black intervals  on their respective sides
(see Fig~\ref{F:equilat}).
This is because if, say $u_1$ is in the blue and $u_2$ is in the black, then  if $u_3$ is in the black, its at distance to $u_1$ is at least $\ell/2$,
 otherwise if  $u_3$ is in the blue then its at distance  to $u_2$ is at least $\ell/2$, in either case contradicting the assumption \eqref{E:contr}.

 Let us therefore assume that all three points are in the blue intervals.
 For simplicity denote $a=\ell/2$,  and by $x_1,x_2,x_3$ the  distances from $u_1,u_2,u_3$  to the midpoints of their respective sides.
 We claim that $d(u_1,u_2)< \frac \ell 2$ implies $x_2<x_1$. Indeed, the law of cosines,
 applied to the triangle with sides $a-x_1$, $a+x_2$, angle of size $\frac{\pi}3$ between them,
 and the remaining side of length $d(u_1,u_2)$, leads to
 \[(a-x_1)^2+(a+x_2)^2-(a-x_1)(a+x_2)<a^2.   \]
 Therefore
 \[ x_1^2+x_2^2+x_1x_2 +ax_2<ax_1.      \]
 Dropping the positive terms $x_1^2+x_2^2+x_1x_2 $ this leaves $ax_2<ax_1$, i.e. $x_2<x_1$.

 Similarly we obtain $x_3<x_2$ and $x_1<x_3$. Since last three inequalities cannot hold simultaneously,
 it means that our assumption \eqref{E:contr} is false.
 \endproof

\remove{\begin{remark} The values given in
Proposition~\ref{P:rectangular} and Proposition~\ref{P:triangular}
estimate from below the support numbers of the rectangular grid
$P_m\times P_n$ and triangular grid $T_n$. We  believe that these
estimates are actually sharp, i.e. these values represent the
largest $k$'s for which the
 respective grids contain $k$-supported cycles.
\end{remark}}

\begin{proof}[Proof of Theorem~\ref{T:kSuppImplStr}]
Let $C$ be a $k$-supported cycle and let $T$  be an arbitrary spanning tree in $G$.
We view the union of $T$ and $C$ as a metric space consisting of unit length intervals, one interval for each edge.
We denote this space by $U$.
It is important to keep in mind that although we use the same notation,
when the embedding of $C$ in $G$ is considered, the cycle $C$ is viewed only as a subgraph, while in the metric space $U$,
$C$ is viewed as a curve (topologically a circle).

In $U$, consider the mapping that takes each edge of $C$ continuously onto its detour in $T$, i.e. the unique path in the
tree $T$ joining the ends of the edge.% (by the definition of the spanning tree such a path exists and is unique).
We define this mapping so that points on an edge of $C$ are mapped
on the detour of that edge in such a way that the ratio of the
distances from a point $t$ inside an edge to the endpoints of the
edge and from its image $t'$ in the detour of the edge to the ends
of the detour path are equal.

If the edge is already in $T$, it is therefore mapped onto itself
identically. Combining these maps we get a continuous map of $C$
into $T$. We emphasize that this happens in the metric space $U$.

Since $C$ is a $k$-supported cycle, it can be partitioned into
paths $I_1, I_2$, and $I_3$ satisfying the conditions of
Definition \ref{D:kSuppCyc}. We color the intervals $I_1, I_2$,
and $I_3$ into blue, green, and red, respectively. To each vertex
that is common to two of the intervals we associate two colors,
one from each interval. Except for these three points, any other
point on $C$ has a unique color associated to it.

\ctikzfig{fig1}

To the continuous map described above from $C$ (thus colored) to
$T$ we apply the following proposition.
\begin{prop}[\cite{RR98} Proposition~5.2]
Let $f:\mathbb{S}^1 \to T$ be a continuous map, and let $\{I_1, I_2, I_3\}$
be an arbitrary partition of $\mathbb{S}^1$ into three intervals with mutually disjoint interiors.
Then there exists $c\in T$ such that $f^{-1}(c)$ has a representative in each of these intervals.
\end{prop}

Let $x,y,z$ be three points, one in each of the $I_1, I_2, I_3$,
mapped onto the same point $c$ in $T$. We note that it is possible
that two of the points coincide, for example $x=y$ if they are at
the intersection of $I_1$ and $I_2$.

If  $x$ is on a blue edge of $C$, and this edge is also in $T$,
then $x=c$. Otherwise, if $x$ is on a blue edge $e$ of $C$ and $e$
is not in $T$, then the cycle formed by $e$ and its detour in $T$
has length at most $\str(G:T)+1$.

In either case, denote by $u_{1}$ the endpoint of $e$ closest to $c$ (select arbitrarily in the case the
endpoints of $e$ are equidistant to $c$ in $U$).
Therefore, $u_{1}$ is a vertex of $C$ colored in blue and  $d_U(u_{1},c)\le \str(G:T)/2$.

Similarly, we have $d_U(u_{2},c)\le  \str(G:T)/2$ and $d_U(u_{3},c)\le  \str(G:T)/2$ for vertices
$u_{2}$ and $u_{3}$ colored in green and red, respectively.
Therefore by triangle inequality each of the distances satisfies
\[ d_C(u_{i},u_{j}) = d_U(u_{i},u_{j})  \le  d_U(u_{i},c) + d_{U}(u_{j},c) \le \str(G:T).\]
Since \[k\le \max_{i,j} d_G(u_{i},u_{j}) \le \max_{i,j} d_C(u_{i},u_{j}) \le \str(G:T),\]
and $T$ was arbitrary, the theorem follows.
\end{proof}

\section{Isometric cycles}
\label{S:IsomCyc}

A somewhat surprising result is obtained by Lokshtanov
\cite[Theorem 3.8]{Lok09}: the length of the longest
 isometric cycle in a finite simple graph can be found in polynomial time.
 From our perspective, the result is interesting because of Corollary~\ref{C:iso_embd}
 which connects the NP-hard notion of stretch with that of polynomially computable maximal length of isometrically embedded cycles.
 We now present a part of  Lokshtanov's construction and make a correction in one of the steps in his argument.

In Lokshtanov's construction, to a graph $G= (V,E)$ and to an integer $k\ge 3$ one
 associates an auxiliary graph $G_k$ as follows. The set of vertices of $G_k$ is the set of ordered pairs
\[ \left\{(u,v) \in V\times V : d(u,v) =
 \left\lfloor \frac k2 \right\rfloor \right\}.\]
 There is an edge in $G_k$ between vertices $(u,v)$ and $(w,x)$ if and only if $(u,w)\in E$ and $(v,x) \in E$.

 For a vertex $(u,v) \in V(G_k)$, define the set
 $$M_k(u,v)=\{ (u,v)\}\ \ \  \mbox{ if $k$ is even, } $$
  and
 $$M_k(u,v)=\{ (u,x)\in V\times V : (u,x) \in V(G_k)  \mbox{ and } \ (v,x) \in E\} \ \ \
 \mbox{ if $k$ is odd.}$$
Lokshtanov  stated the following theorem.
\begin{theorem}[\cite{Lok09}]
\label{T:IsomC}
A graph $G$ has an isometric cycle of length $k$ if and only if there are vertices $u$ and $v$ and $x$ in $V(G)$
so that $(v, x) \in  M_k(v, u)$ and
$d_{G_k} [(u, v), (v, x)] = \left\lfloor \frac k2 \right\rfloor$.
\end{theorem}
The theorem is proved in two cases, for $k$ even and for $k$ odd.
For odd $k$ the proof is based on the following lemma.

\begin{lemma}[{\cite[Lemma 3.6]{Lok09}}]\label{L:Lok} If $k$ is odd and the distance in $G_k$ between $(u,v)$ and a vertex $(v,x)$ in
$M_k(v,u)$ is $\left\lfloor \frac k2 \right\rfloor$, then there is
an isometric cycle of length $k$ in $G$, going through $u$ and
$v$.
\end{lemma}

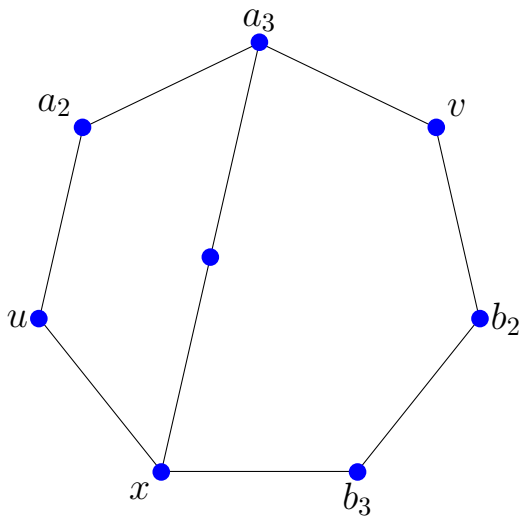
\begin{figure}[h!]
    \centering

    \begin{tikzpicture}
    \node[draw, regular polygon, regular polygon sides=7, minimum size=6cm] (a) {};
     \draw (a.corner 4)--(a.corner 1);
    \foreach \i in {1,2,...,7}
        \node[fill=blue, circle, scale=0.6] at (a.corner \i){};
    \node[above] at (a.corner 1){\large $a_3$};
    \node[above right] at (a.corner 7){\large $v$};
    \node[right] at (a.corner 6){\large $b_2$};
    \node[below] at (a.corner 5){\large $b_3$};
    \node[below left] at (a.corner 4){\large $x$};
    \node[left] at (a.corner 3){\large $u$};
    \node[above left] at (a.corner 2){\large $a_2$};
    \node[fill=blue, circle, scale=0.6] at ($(a.corner 4)!0.5!(a.corner 1)$) {};
    \end{tikzpicture}
     \caption{Counterexample to Lemma 3.6 in [Lok09]}\label{fig:ceL}
\end{figure}

We should remark that Lemma~3.6 in \cite{Lok09} is erroneous, as
the graph in Figure \ref{fig:ceL} shows. In fact, the graph has no
isometric cycle of length $7$, while all conditions of the Lemma
\ref{L:Lok} are satisfied.

The case of  even $k$  in Theorem \ref{T:IsomC} (whose proof in
\cite{Lok09} is correct) is contained in
\cite[Corollary~3.4]{Lok09}. We restate it  as
\begin{theorem}[\cite{Lok09}]
If $k$ is even, there is an isometric cycle of length $k$ in $G$ if and only if there is a pair of vertices $u$ and $v$ with
$d_{G_k} [(u, v), (v, u)] =  k/2$.
\label{T:kev}
\end{theorem}

Our goal now is to prove:

\begin{observation}\label{O:LokCorr} Lokshtanov's theorem for odd $k$ can be
derived from Theorem \ref{T:kev}.
\end{observation}

\begin{proof} We show that the existence of an isometric cycle of odd
length $k$ in $G$ is equivalent to the existence of an isometric
cycle of length $2k$ is an auxiliary bipartite graph $G'$.

We derive the auxiliary graph $G'$ from a graph $G$ as follows.
Denote $|V(G)|=n$ and $|E(G)|=m$ the number of vertices and edges
in the connected simple graph $G$. For visualization purposes
consider the vertices of $G$ colored in blue. Define $G'$ from $G$
by dividing every edge in $E(G)$ into two edges in $E(G')$ by the
introduction of a new red vertex. The set $V(G')$ thus consists of
$n+m$ vertices, $n$ of which are blue and $m$ are red. There are
$2m$ edges in $E(G')$, all joining a blue and a red vertex,
therefore $G'$ is bipartite. Note that for any two blue vertices
$u$ and $v$, we have
\[ d_{G'}(u,v) =2 d_{G}(u,v).\]
We use Theorem~\ref{T:kev} to find isometric cycles of length $2k$
in the graph $G'$. Since in such a cycle $C'_{2k}$ blue and red
vertices alternate, we can label the vertices of the cycle as
\[b_0, r_1, b_1, r_2, b_2, \dots, b_{k-1}, r_{k}, b_{k}=b_0.\] By
dropping the red vertices, we obtain the cycle $C_k$ in $G$ with
vertices $b_0, b_1, b_2, \dots, b_k=b_0$. If for two vertices $u,
v$ in $C_k$ we had
 \[ d_{G}(u,v) < d_{C_k}(u,v),\]
 this would imply
 \[ d_{G'}(u,v) < 2 d_{C_k}(u,v) = d_{C'_{2k}}(u,v).\]
 But this contradicts the fact that $C'_{2k}$ is isometric in $G'$.
 Therefore the cycle $C_k$ embeds isometrically in $G$.\end{proof}

Our next goal is the following observation.

\begin{observation} In some cases
Corollary~\ref{C:iso_embd} provides very weak estimates for
stretch.
\end{observation}

In fact, by Theorem \ref{T:kSuppImplStr}, Proposition
\ref{P:rectangular} implies that $\sigma(P_n\times P_n)\ge n-1$
and Proposition \ref{P:triangular} implies that
$\sigma(T_n)\ge\lfloor\frac{n-1}2\rfloor$. On the other hand one
can show that the maximal sizes of isometric cycles in $P_n\times
P_n$ and $T_n$ are $4$ and $3$ respectively.

For $P_n\times P_n$ it is proved as follows:
Clearly, the cycle of length four is isometrically embedded. This is the only case when two consecutive turns of either
 $90^\circ$  or  $-90^\circ$  are allowed for edges $x_0x_1$, $x_1x_2$, and $x_2x_3$.
Assume $x_0,x_1,\dots,x_k=x_0$ is an isometrically embedded cycle of length $k>4$.
For any $3 \le i\le k$, if there are turns in the path between $x_0$ and $x_i$, they have to alternate so that after a turn by
$90^\circ$ the next turn will be by $-90^\circ$ and vice versa.
However, such alternating turns only force the distance $d(x_0,x_i)$ to increase and therefore the path cannot close into a cycle.

For $T_n$ it is proved as follows:
Clearly, the cycle of length three is isometrically embedded. This is the only case when a turn of either
 $60^\circ$  or  $-60^\circ$  is allowed for edges $x_0x_1$, $x_1x_2$.
Assume $x_0,x_1,\dots,x_k=x_0$ is an isometrically embedded cycle of length $k>3$.
For any $2 \le i\le k$, if there are turns in the path between $x_0$ and $x_i$, they have to alternate so that after a turn by
$120^\circ$ the next turn will be by $-120^\circ$ and vice versa.
However, such alternating turns only force the distance $d(x_0,x_i)$ to increase and therefore the path cannot close into a cycle.

\section{Approximation algorithm for $k$-supported cycles}
\label{S:approx_k_sp_Cy}

Our goal in this section is to construct an approximation
algorithm (in the sense of \cite{Vaz01}, \cite{WS11}) for computing the support number, i.e. the
maximal $k$ such that a given graph $G$ contains a $k$-supported
cycle.

First we define a polynomial algorithm for computing the integer
characteristic $W(G)$ of a graph $G$ introduced in
Definition~\ref{D:cw} below. Then we prove Theorems~\ref{T:(A)}
and \ref{T:(B)}.

For each vertex $r$ in a simple connected graph $G=(V(G),E(G))$ we build the breadth-first
search tree (see \cite{BM08} pp.~137--139) rooted at $r$.
For $x\in V(G)$ denote the {\it  level} of $x$ by $\ell(x) = d_G(r,x)$.
For each $n\in \mathbb{N}$ we consider the  subset  $R(n)$ of edges in $G$
 consisting of edges $xy$ for which
    \[\max\{\ell(x),\ell(y)\}> n.\]
Thus, $R(n)$ is the set of edges of $G$ with at least one endpoint
outside the ball centered at the root $r$ of radius $n$.

\begin{definition}\label{D:cw}  Write $x\simeq_n y$ is $x$ and  $y$ in the same component of
$(V(G),R(n))$.  Define  the \emph{cycle width} of $G$ as $W(G)=
\displaystyle{\max_{r\in V(G)} W(r)}$ where
\[W(r)= \max_n\max\{d_G(x,y):~ \ell(x)=\ell(y)=n,~ x\simeq_n y\}.\]
\end{definition}

\begin{remark}
Since the distance matrix for $V(G)$ is polynomially computable, it is straightforward that  $W(G)$ is also a polynomially computable
quantity.
\end{remark}

\begin{proof}[Proof of Theorem \ref{T:(A)}]
Suppose that we have computed $W(G)$ and found a root $r$, a radius $n$, and
corresponding vertices  $x$ and $y$ such that
\[ \ell(x)=\ell(y)=n, \mbox{ and }     d_{G}(x,y) = W(G).\]
Denote by $T$ the breadth-first tree rooted at $r$.
To produce a $k$-supported cycle with $k\ge W(G)/3$, we start at the point $w$ where the paths $rx$ and $ry$ in $T$
depart from each other (it can be the root $r$), and consider the cycle
consisting of
\begin{itemize}
\item path from $x$ to $w$
\item path from $w$ to $y$
\item path consisting of edges in $R(n)$ joining $x$ and $y$.
\end{itemize}
The paths above form the partition $I_1$, $I_2$, and $I_3$ of the cycle under consideration.
Consider arbitrary vertices $u_{1}\in I_{1}$, $u_{2}\in I_{2}$, and $u_{3}\in I_{3}$.
Because $\ell(u_{3}) \ge n$, we have that
\[ d_{G}(u_1,u_3) \ge d_{G}(u_1,x) \ \mbox{ and }  \  d_{G}(u_2,u_3) \ge d_{G}(u_2,y).\]
Since
\[ d_{G}(x, u_1) + d_{G}(u_1,u_2) + d_{G}(u_2,y) \ge  d_{G}(x,y) = W(G),\]
we obtain that
\[ d_{G}(u_1,u_3)+ d_{G}(u_1,u_2) +d_{G}(u_2,u_3) \ge W(G).\]
Therefore the maximum of the three distances has to be at least $W(G)/3$ and thus the cycle is $W(G)/3$-supported.
\end{proof}

We now bound the support number of $G$ from above in terms of
$W(G)$.

\begin{proof}[Proof of Theorem \ref{T:(B)}] We need to show that the existence of a $k$-supported cycle in
$G$ implies $W(G)\ge k-4$.

Let $C$ be a $k$-supported cycle in $G$, and let $x$, $y$, $z$ be
the vertices where $C$ is partitioned into $I_1$, $I_2$, and $I_3$
showing the situation. For an arbitrarily selected root $r$ we
prove that the presence of the cycle $C$ implies $W(r)+4 \ge k$
(in the estimates we will use however the larger $W(G)$ instead of
$W(r)$).

We introduce additional structure on the vertex set of $G$, relative to the chosen root $r$.

\begin{enumerate}[{\bf (a)}]

\item We split the vertex set of $G$ into the set of
\emph{levels}, that is, for each $n\in \mathbb{N}$, the
\emph{level $n$} is the set of vertices $u$ satisfying
$d_G(u,r)=n$.

\item For each $n\in\mathbb{N}$ we split the level $n$ into subsets
which are equivalence classes of the following relation: $x\sim y$ if
$x=y$ or if $x$ and $y$ are joined by a path in $R(n)$.
In such a way we get pairwise disjoint sets of vertices. Let us denote them
$\{C(n,j)\}_{j=1}^{c(n)}$ and call them \emph{blocks of level
$n$}.

\item We also introduce \emph{bunches of level $n$} as unions of
blocks of level $n$ defined by the following condition: If there
is an edge joining two blocks of level $n$, they should be in the same bunch.
More formally: Blocks $C(n,j)$ and $C(n,i)$ are in the same bunch
if and only if there exists a sequence of blocks
$\{C(n,j(t))\}_{t=0}^{N(j,i)}$ such that $C(n,j(0))=C(n,j)$,
$C(n,j(N(j,i)))=C(n,i)$ and there is an edge between $C(n,j(t))$
and $C(n,j(t+1))$ for each $t = 0, \dots, N(j,i)-1$.

We denote bunches of level $n$ by $\{B(n,s)\}_{s=1}^{b(n)}$.

\end{enumerate}

\begin{lemma} \label{L:btb}
All edges of $G$ going from a bunch $B(n,s)$ to the level $(n-1)$ have
their other ends in the same block of level $(n-1)$.
\end{lemma}

\begin{proof} If $u$ and $v$ are such ends in level $(n-1)$ we can
construct a path in $R(n-1)$ joining them.
\end{proof}

\begin{corollary}\label{C:W(G)+2}
The diameter of a bunch is at most $W(G)+2$.
\end{corollary}
\begin{proof}
By the definition of cycle width $W(G)$, any block of any level has diameter at most $W(G)$.
Since by Lemma~\ref{L:btb} the  predecessors  of a bunch at level $n$ are all in  a block at level $(n-1)$,
we obtain that the diameter of the bunch is no more than $W(G)+2$.
\end{proof}

Given a bunch or a block, by the part of $G$ ``above'' the bunch
or block we mean the set of vertices of $G$ which are disconnected
from the root by the removal of the bunch or block.

Assume that the graph contains a $k$-supported cycle $C$, and
color vertices of $C$ blue, green, and red according to the
partition of $C$  into $I_1$, $I_2$ and $I_3$ (following
Definition \ref{D:kSuppCyc}). Vertices where two of the intervals
meet are considered to have both colors.

Corollary \ref{C:W(G)+2} implies that if there is a bunch which
contains points of all three colors, then $k\le W(G)+2$.

We prove that there is a bunch containing vertices of two colors
and having distance $1$ to a vertex of the third color. We start
by proving the following lemma.

\begin{lemma} Either there exists a bunch containing all three
colors, or there exists a bunch such that one of the colors is
present only ``above'' the bunch.
\end{lemma}

\begin{proof}
Consider three cases:

\begin{enumerate}[{\bf (1)}]

\item The cycle contains the root.

It is easy to see that in this case the cycle intersects exactly
one of the bunches of level $1$. If the bunch does not contain all
three colors, then one of the colors is only ``above'' the bunch.

\item The cycle does not contain the root, but contains a vertex
of level $1$.

It is easy to see that in this case the cycle intersects only one
of the bunches of level $1$. Unless the intersection contains all
three colors, one of the colors is only ``above'' the bunch.

\item The lowest-level vertices in the cycle have level $2$ or
higher.

In this case all three colors are ``above'' the block immediately
preceding the cycle (we mean the block of the previous level
through which the cycle is connected with the root).

\end{enumerate}
\end{proof}

To conclude the proof of the Theorem~\ref{T:(B)}, we only have to consider
the case when no bunch contains all three colors.

Let bunch $B(n,s)$ be the bunch with the largest $n$ such that one
of the colors is present only ``above'' the bunch. In such a case
$B(n,s)$ together with those vertices of level $(n+1)$ which are
``above'' $B(n,s)$ should contain all three colors. In fact, since
any two colors should ``meet'', otherwise there are bunches
``above'' $B(n,s)$ for which one of the colors is present only
``above'' the bunch.

Since the diameter of  $B(n,s)$ is at most $W(G)+2$, we conclude
that there is a triple of points of different colors with pairwise
distances $\le W(G)+4$.
\end{proof}

\end{document}